\documentclass[a4paper]{article}
\usepackage{amssymb}
\usepackage{graphicx}
\usepackage{hyperref}
\usepackage{latexsym}
\newtheorem{theorem}{Theorem}
\newtheorem{corollary}{Corollary}
\newtheorem{lemma}{Lemma}
\newtheorem{definition}{Definition}
\newtheorem{thesis}{Thesis}

\begin{document}
\title{How definitive \textit{is} the standard interpretation of Goodstein's argument}
\author{\normalsize{Bhupinder Singh Anand}}
\date{\small{Minor editing; new Lemma \ref{sec:lem6}; graph of terminating Goodstein sequences included; major revision of notation on \today.\footnote{Subject class: LO; MSC: 03B10}}}
\maketitle

\begin{abstract}
\noindent Goodstein's argument is essentially that the hereditary representation $m_{[b]}$ of any given natural number $m$ in the natural number base $b$ can be mirrored in Cantor Arithmetic, and used to well-define a finite decreasing sequence of transfinite ordinals each of which is not smaller than the ordinal corresponding to the corresponding member of Goodstein's sequence of natural numbers $G(m)$. The standard interpretation of this argument is first that $G(m)$ must therefore converge; and second that this conclusion---Goodstein's Theorem---is unprovable in Peano Arithmetic but true under the standard interpretation of the Arithmetic. We argue however that even assuming Goodstein's Theorem is indeed unprovable in PA, its truth must nevertheless be an intuitionistically unobjectionable consequence of some constructive interpretation of Goodstein's reasoning. We consider such an interpretation, and construct a Goodstein functional sequence to highlight why the standard interpretation of Goodstein's argument ought not to be accepted as a definitive property of the natural numbers.

\vspace{+1.5ex}
\end{abstract}

\section{Introduction}
\label{intro:1}

Goodstein's argument \cite{Gd44} is essentially that the hereditary representation $m_{[b]}$ in Peano Arithmetic\footnote{By Peano Arithmetic we mean the arithmetic of the natural numbers that is definable in a formal number theory such as Mendelson's theory S (\cite{Me64}, p102).}, of any given natural number $m$ in the natural number base $b$ can be mirrored in Cantor's (ordinal) Arithmetic\footnote{By Cantor Arithmetic we mean the arithmetic of the ordinal numbers that is definable in a formal set theory such as ZFC or NBG (cf. \cite{Me64}, p189).}, and used to yield a finite\footnote{We note that for a sequence of ordinals to be termed as finite it must be a well-defined set in Cantor Arithmetic (\cite{Me64}, p184).}, decreasing, sequence of transfinite ordinals, each of which is not smaller\footnote{In the sense in which this relation is defined in Cantor Arithmetic.} than the ordinal corresponding to the corresponding member of Goodstein's sequence of natural numbers $G(m)$. 

The standard interpretation of this argument is first that $G(m)$ must therefore converge (Goodstein's Theorem); second that this number-theoretic proposition is unprovable in any formal system of Peano Arithmetic, but expresses a truth under the standard interpretation of the Arithmetic that appeals necessarily to transfinite reasoning (Kirby-Paris Theorem \cite{KP82}); and third that Goodstein's Theorem is, in a sense, a proposition that under such interpretation expresses a more natural independence phenomenon than G\"{o}del's Theorem on formally unprovable, but interpretatively true, sentences of any formal system of Peano Arithmetic.

However we note first that G\"{o}del's reasoning can be carried out in a weak Arithmetic such as Robinson's system Q \cite{Ro50}, which does not admit mathematical induction. The truth of the unprovable G\"{o}del sentence could thus be reasonably argued as being even more intuitive than the truth---under the standard interpretation---of any number-theoretic assertion of Peano Arithmetic that necessarily appeals to mathematical induction. Moreover both truths are classically accepted as constructive and intuitionistically unobjectionable. 

We note further that Goodstein's Theorem involves a non-constructive---hence non-verifiable---concept of mathematical truth\footnote{Necessarily so, according to a reasonable interpretation of the Kirby-Paris Theorem \cite{KP82}.} that, prima facie, is of a higher order of intuition---in a manner of speaking---than that required to see that G\"{o}del's formally unprovable sentence is a true number-theoretic assertion of Peano Arithmetic under its standard interpretation \cite{Go31a}.\footnote{This follows since, as Gödel pertinently notes in his seminal 1931 paper (\cite{Go31a}, p26), the truth of his unprovable sentence, under its standard interpretation, is meta-mathematically verifiable constructively, in an intuitionistically unobjectionable manner, in Peano Arithmetic. See, also, the reasoning in \cite{An02a}.}  

If, therefore, the proof of Goodstein's Theorem is to be considered as having established both an unprovable proposition of Peano Arithmetic that is true under its standard interpretation, and a more natural independence phenomenon than G\"{o}del's, then such truth too must reasonably be a consequence of some constructive---and intuitionistically unobjectionable---interpretation of Goodstein's reasoning. 

In Section \S\ref{thesis:1} we argue that such an interpretation does indeed exist as an implicit thesis in Goodstein's argument. In Sections \S\ref{func:1} and \S\ref{obj:1} we consequently construct and consider a Goodstein functional sequence that highlights why the standard interpretations of this argument ought not to be considered as definitive.

\section{Notation and Definitions}
\label{def:1}

\noindent \textbf{Ordinal number notation}: We denote the ordinal number corresponding to the natural number $m$ by $[m]$.

\vspace{+1ex}
\noindent \textbf{Hereditary representation $m_{[b]}$ of a natural number $m$ in base $b$}: The hereditary representation of the natural number $m$ in the natural number base $b$\footnote{It is implicit here, and in what follows, that (the base) $b > 1$.}, which we denote by $m_{[b]}$\footnote{We note that $m$ and $m_{[b]}$ denote, but are different syntactical expressions of, the same natural number.}, is its syntactic expression as a sum of powers of the natural number base $b$, followed by expression of each of the exponents as a sum of powers of $b$ etc., until the process stops.

\vspace{+1ex}
\noindent \textbf{The rank of a hereditary representation}: The rank of a hereditary representation is the highest power of the natural number base that has a non-zero co-efficient in the representation.

\vspace{+1ex}
\noindent \textbf{Goodstein Sequence}: Let $m_{[b\ \hookrightarrow\ (b+1)]}$ be the non-negative integer which results if we syntactically replace each $b$ by $(b+1)$ in the hereditary representation $m_{[b]}$ of a natural number $m$ in base $b$\footnote{We note that $m_{[b\ \hookrightarrow\ (b+1)]}$ is also, sometimes, denoted in the literature by $B(b, m)$, where $B$ is referred to as the Goodstein `bumping' or `base-change' function, and $B(b, m)$ is the non-negative integer that results if we syntactically replace each $b$ by $(b+1)$ in the hereditary representation $m_{[b]}$ of the natural number $m$ in the base $b$.}. Starting at the hereditary representation of the natural number $m$ in base $2$, we formally define the Goodstein sequence, $G(m)$ as:

\begin{eqnarray}
\{g_{1}(m),\ g_{2}(m),\ g_{3}(m) \ldots\}
\end{eqnarray}

\vspace{+1ex}
\noindent where:

\begin{eqnarray}
g_{n}(m) & \equiv & (g_{(n-1)}(m)_{[n\ \hookrightarrow\ (n+1)]} -1)
\end{eqnarray}

\vspace{+1ex}
\noindent \textbf{Goodstein's Theorem}: For all natural numbers $m$, there exists a natural number $n$ such that the $n$th term, $g_{n}(m)_{[(n+1)]}$, of the Goodstein sequence $G(m)$ is 0. 
 
\section{The Goodstein operations}
\label{opera:1}

We note that each natural number $m$ has a unique hereditary representation, of some finite rank $l$, in any given natural number base $b$. Without loss of generality, we may express this syntactically as:

\begin{eqnarray}
m_{[b]} & = & \sum_{i=0}^{l} a_{i}.b^{i_{[b]}}
\end{eqnarray}

\vspace{+1ex}
\noindent where:

\begin{tabbing}
(a) \= $0 \leq a_{i} < b$ over $0 \leq i \leq l$; \\
(b) \= $a_{l} \neq 0$; \\
(c)	\= for each $0 \leq i \leq l$ the exponent $i$ too is expressed syntactically \\ \> by its hereditary representation $i_{[b]}$ in the base $[b]$; as also are all \\ \> of its exponents and, in turn, all of their exponents, etc.
\end{tabbing}

\subsection{Goodstein bumping operation}
\label{opera:1.1}

We define the Goodstein bumping operation\footnote{We prefer to define these concepts as `operations', rather than as `functions', simply to avoid any implicit commitment to possible set-theoretic properties.}, on the hereditary representation of the natural number $m$ in the natural number base $b$, by:

\vspace{+1ex}
\textbf{The Goodstein bumping operation}: The natural number $m_{[b\ \hookrightarrow\ (b+1)]}$ is derived from $m$ under a Goodstein bumping operation by syntactically replacing the base $b$ by the base $(b+1)$ in the hereditary representation $m_{[b]}$ of $m$ in the base $b$. 

\vspace{+1ex}
We can express this, also without loss of generality, as:

\begin{eqnarray}
m_{[b\ \hookrightarrow\ (b+1)]} & = & \sum_{i=0}^{l} a_{i}.(b+1)^{i_{[b\ \hookrightarrow\ (b+1)]}}
\end{eqnarray}

\subsection{Complete Goodstein operation}
\label{opera:1.2}

We, then, define the result of a complete Goodstein operation, on the hereditary representation of a natural number $m$ in the natural number base $b$, as the natural number $m_{[b\ \hookrightarrow\ (b+1)]} -1$. 

\section{Goodstein's argument}
\label{argue:1}

In his 1944 paper Goodstein essentially considers, for any given natural number $m$, the sequence, say $G(m)$, of natural numbers of Peano Arithmetic:

\begin{eqnarray}
\{g_{1}(m)_{[2]},\ g_{2}(m)_{[3]},\ g_{3}(m)_{[4]},\ g_{4}(m)_{[5]},\ \ldots\}
\end{eqnarray}

\vspace{+1ex}
\noindent and the parallel sequence, $O(m_{o})$, of ordinal numbers of Cantor Arithmetic:

\begin{eqnarray}
\{g_{1}(m_{o})_{[2_{o}\ \hookrightarrow\ \omega]},\ g_{2}(m_{o})_{[3_{o}\ \hookrightarrow\ \omega]},\ g_{3}(m_{o})_{[4_{o}\ \hookrightarrow\ \omega]},\ g_{4}(m_{o})_{[5_{o}\ \hookrightarrow\ \omega]},\ \ldots\}
\end{eqnarray}

\vspace{+1ex}
\noindent where $g_{n}(m_{o})_{[(n_{o}+1_{o})\ \hookrightarrow\ \omega]}$ is the ordinal number obtained from $g_{n}(m)_{[(n+1)]}$ by syntactically replacing all natural numbers in $g_{n}(m)_{[(n+1)]}$ by their corresponding ordinals, other terms by their corresponding\footnote{We assume that such a correspondence exists, along the lines described, for instance, in Mendelson (\cite{Me64}, p192).} set-theoretical terms, and then syntactically replacing the ordinal base $[n_{o}+1_{o}]$ by the ordinal $[\omega]$. 

\vspace{+1ex}
Now, by properties of ordinal addition, multiplication and exponentiation (\cite{Me64}, p189) we have, for any given natural number $m_{[b]}$, the ordinal inequality, $m_{o}{_{[b_{o}\ \hookrightarrow\ \omega]} > m_{o}{_{[b_{o}]}}}$, where $m_{o}$ denotes the ordinal corresponding to the natural number $m$, and $\omega$ denotes Cantor's lowest transfinite ordinal. 

Further, by arithmetical properties that are characteristic of transfinite ordinals such as $\omega$, it can be shown that $O(m_{o})$ is a decreasing sequence of transfinite ordinals, each of whose members is not smaller than the ordinal corresponding to the corresponding member of the sequence of natural numbers, $G(m)$.

Since in Cantor Arithmetic the ordinals are well-ordered, and there are no infinite, decreasing, sequences of ordinals, Goodstein concludes that $O(m_{o})$ is a finite sequence of transfinite ordinals. 

\section{Goodstein's Theorem}
\label{theo:1}

Now, assuming that a set theory such as ZFC can be treated as a consistent extension of the first order Peano Arithmetic PA\footnote{That a set theory, such as ZFC, cannot be unrestrictedly treated as a consistent extension of a first order Peano Arithmetic, such as standard PA, is also suggested by independent arguments offered in \cite{An02b}.}, the standard interpretation of Goodstein's argument is, then: 

\vspace{+1ex}
\textbf{Goodstein's Theorem}: For all natural numbers $m$, there exists a natural number $n$ such that the $n^{th}$ term $g_{n}(m)$ of the Goodstein sequence $G(m)$ is 0.

\vspace{+1ex}
We note that this interpretation of Goodstein's argument is supported by the following examples.

\subsection{Example: $m = 1$}
\label{theo:1.1}

\begin{quote}
$g_{1}(1)_{[2]} = 1 \cdot 2^{0}$

\begin{quote}
\footnotesize $g_{1}(1)_{[2\ \hookrightarrow\ 3]} = 1 \cdot 3^{0}$
\end{quote}

$g_{2}(1)_{[3]} = 0 \cdot 3^{0}$
\end{quote}

\vspace{+1ex}
Hence $G(1)$ is $\{1,\ 0\}$.

\subsection{Example: $m = 2$}
\label{theo:1.2}

\begin{quote}
$g_{1}(2)_{[2]} = 1\cdot2^{1} + 0\cdot2^{0}$

\begin{quote}
\footnotesize $g_{1}(2)_{[2\ \hookrightarrow\ 3]} = 1\cdot3^{1} + 0\cdot3^{0}$
\end{quote}

$g_{2}(2)_{[3]} = 0\cdot3^{1} + 2\cdot3^{0}$

\begin{quote}
\footnotesize $g_{2}(2)_{[3\ \hookrightarrow\ 4]} = 0\cdot4^{1} + 2\cdot4^{0}$
\end{quote}

$g_{3}(2)_{[4]} = 0\cdot4^{1} + 1\cdot4^{0}$

\begin{quote}
\footnotesize $g_{3}(2)_{[4\ \hookrightarrow\ 5]} = 0\cdot5^{1} + 1\cdot5^{0}$
\end{quote}

$g_{4}(2)_{[5]} = 0\cdot5^{1} + 0\cdot5^{0}$
\end{quote}

\vspace{+1ex}
Hence $G(2)$ is $\{2,\ 2,\ 1,\ 0\}$.
 
\subsection{Example: $m = 3$}
\label{theo:1.3}

\begin{quote}
$g_{1}(3)_{[2]} = 1\cdot2^{1} + 1\cdot2^{0}$

\begin{quote}
\footnotesize $g_{1}(3)_{[2\ \hookrightarrow\ 3]} = 1\cdot3^{1} + 1\cdot3^{0}$
\end{quote}

$g_{2}(3)_{[3]} = 1\cdot3^{1} + 0\cdot3^{0}$

\begin{quote}
\footnotesize $g_{2}(3)_{[3\ \hookrightarrow\ 4]} = 1\cdot4^{1} + 0\cdot4^{0}$
\end{quote}

$g_{3}(3)_{[4]} = 0\cdot4^{1} + 3\cdot4^{0}$

\begin{quote}
\footnotesize $g_{3}(3)_{[4\ \hookrightarrow\ 5]} = 0\cdot5^{1} + 3\cdot5^{0}$
\end{quote}

$g_{4}(3)_{[5]} = 0\cdot5^{1} + 2\cdot5^{0}$

\begin{quote}
\footnotesize $g_{4}(3)_{[5\ \hookrightarrow\ 6]} = 0\cdot6^{1} + 2\cdot6^{0}$
\end{quote}

$g_{5}(3)_{[6]} = 0\cdot6^{1} + 1\cdot6^{0}$

\begin{quote}
\footnotesize $g_{5}(3)_{[6\ \hookrightarrow\ 7]} = 0\cdot7^{1} + 1\cdot7^{0}$
\end{quote}

$g_{6}(3)_{[7]} = 0\cdot7^{1} + 0\cdot7^{0}$
\end{quote}

\vspace{+1ex}
Hence $G(3)$ is $\{3,\ 3,\ 3,\ 2,\ 1,\ 0\}$.

\section{Some Goodstein sequence lemmas}
\label{lem:1}

However, we now argue that to the extent that standard interpretations of Goodstein's argument use, but ignore the significance of, the fact that arithmetical properties of the natural number sequence $G(m)$ are not necessarily shared by the corresponding transfinite ordinal sequence $O(m_{o})$, such interpretations cannot be considered definitive. 

In order to highlight the significance of the above distinction, we introduce some general properties of sequences generated by iterated application of the complete Goodstein operation on the hereditary representation of the natural number $m$ in a natural number base $b$.

\subsection{The first Goodstein sequence lemma}
\label{lem:1.1}

We note first that:

\begin{lemma}
\label{sec:lem1}
If there exists a natural number $n$ such that the $n^{th}$ term $g_{n}(m)$ of the Goodstein sequence $G(m)$ is $0$ then, for all $m > 3$, there is some $k < n$ such that $g_{(k-1)}(m) > g_{(k)}(m)$.
\end{lemma}

\vspace{+1ex}
\noindent \textit{Proof}: We have by definition:

\begin{eqnarray}
g_{(n-1)}(m)_{[n]} & \equiv & \sum_{i=0}^{l} a_{i}.n^{i_{[n]}}
\end{eqnarray}

\vspace{+1ex}
\noindent and:

\begin{eqnarray}
g_{n}(m)_{[(n+1)]} & \equiv & (\sum_{i=0}^{l} a_{i}.(n+1)^{i_{[n\ \hookrightarrow\ (n+1)]}}) - 1
\end{eqnarray}

\vspace{+1ex}
\noindent where:

\begin{tabbing}
(a) \= $0 \leq a_{i} < n$ over $0 \leq i \leq l$; \\
(b) \= $a_{l} \neq 0$; \\
(c)	\= for each $0 \leq i \leq l$ the exponent $i$ too is expressed syntactically \\ \> by its hereditary representation $i_{[n]}$ in the base $[n]$; as also are all \\ \> of its exponents and, in turn, all of their exponents, etc.
\end{tabbing}

\noindent Hence:

\begin{eqnarray*}
g_{2}(m)_{[3]} - g_{1}(m)_{[2]} & \equiv & \sum_{i=0}^{l} (a_{i}.3^{i_{[2\ \hookrightarrow\ 3]}} - a_{i}.2^{i_{[2]}}) - 1
\end{eqnarray*}

\vspace{+1ex}
\noindent Now, if $m > 3$ then $l \geq 2$. Hence:

\begin{eqnarray*}
g_{2}(m)_{[3]} - g_{1}(m)_{[2]} & \geq & (3^{l_{[2\ \hookrightarrow\ 3]}} - 2^{l_{[2]}}) - 1 \\
& \geq & (3^{2_{[2\ \hookrightarrow\ 3]}} - 2^{2_{[2]}}) - 1 \\
& > & 1
\end{eqnarray*}

\vspace{+1ex}
\noindent The lemma follows. \hfill $\Box$

\subsection{The second Goodstein sequence lemma}
\label{lem:1.2}

\begin{lemma}
\label{sec:lem2}
If the hereditary representation of the $k^{th}$ term, $g_{k}(m)$ of the Goodstein sequence $G(m)$ contains more than one non-zero term, then $g_{(k+1)}(m) \geq g_{k}(m)$.
\end{lemma}

\vspace{+1ex}
\noindent \textit{Proof}: We have:

\begin{eqnarray*}
g_{(k+1)}(m)_{[(k+2)]} - g_{k}(m)_{[(k+1)]} & \equiv & (\sum_{i=0}^{l} a_{i}((k+2)^{i_{[(k+1)\ \hookrightarrow\ (k+2)]}} - (k+1)^{i_{[(k+1)]}})) - 1
\end{eqnarray*}

\vspace{+1ex}
\noindent where:

\begin{tabbing}
(a) \= $0 \leq a_{i} < k$ over $0 \leq i \leq l$; \\
(b) \= $a_{l} \neq 0$; \\
(c)	\= for each $0 \leq i \leq l$ the exponent $i$ too is expressed syntactically \\ \> by its hereditary representation $i_{[k]}$ in the base $[k]$; as also are all \\ \> of its exponents and, in turn, all of their exponents, etc.
\end{tabbing}

If now $0 \leq j < l$ and $a_{j}, a_{l} \neq 0$, then $l \geq 1$ and:

\begin{eqnarray*}
g_{(k+1)}(m)_{[(k+2)]} - g_{k}(m)_{[(k+1)]} & \geq & a_{l}((k+2)^{l_{[(k+1)\ \hookrightarrow\ (k+2)]}} - (k+1)^{l_{[(k+1)]}}) - 1 \\
& \geq & a_{l}((k+2)^{1_{[(k+1)\ \hookrightarrow\ (k+2)]}} - (k+1)^{1_{[(k+1)]}}) - 1 \\
& \geq & 0
\end{eqnarray*}

\vspace{+1ex}
The lemma follows. \hfill $\Box$

\subsection{The third Goodstein sequence lemma}
\label{lem:1.3}

\begin{lemma}
\label{sec:lem3}
If the hereditary representation of the $k^{th}$ term $g_{k}(m)$ of the Goodstein sequence $G(m)$ is of the form $a_{l}\cdot(k+1)^{l_{[(k+1)]}}$, where $1 \leq a_{l} < (k+1)$ then $g_{(k+1)}(m) \geq g_{k}(m)$ if $l \geq 1$ and $g_{(k+1)}(m) < g_{k}(m)$ if $l = 0$.
\end{lemma}

\vspace{+1ex}
\noindent \textit{Proof}: We have:

\begin{eqnarray*}
g_{k}(m)_{[(k+1)]} & = & a_{l}\cdot(k+1)^{l_{[(k+1)]}}
\end{eqnarray*}

\vspace{+1ex}
\noindent and:

\begin{eqnarray*}
g_{(k+1)}(m)_{[(k+2)]} - g_{k}(m)_{[(k+1)]} & = & a_{l}((k+2)^{l_{[(k+1)\ \hookrightarrow\ (k+2)]}} - (k+1)^{l_{[(k+1)]}}) - 1
\end{eqnarray*}

\vspace{+1ex}
Then if $l \geq 1$:

\begin{eqnarray*}
g_{(k+1)}(m)_{[(k+2)]} - g_{k}(m)_{[(k+1)]} & \geq & a_{l}((k+2) - (k+1)) - 1 \\
& \geq & a_{l} - 1 \\
& \geq & 0
\end{eqnarray*}

\vspace{+1ex}
Whilst if $l = 0$:

\begin{eqnarray*}
g_{(k+1)}(m)_{[(k+2)]} - g_{k}(m)_{[(k+1)]} & = & a_{l}((k+2)^{0} - (k+1)^{0}) - 1 \\
& = & - 1
\end{eqnarray*}

\vspace{+1ex}
The lemma follows. \hfill $\Box$

\subsection{The fourth Goodstein sequence lemma}
\label{lem:1.4}

\begin{lemma}
\label{sec:lem4} If the $n^{th}$ term $g_{n}(m)$ of the Goodstein sequence $G(m)$ is $0$, then $n = 2n_{1}$ for some natural number $n_{1} > 1$ and $g_{(n_{1}-1)}(m)_{[n_{1}]} = 1\cdot n_{1}^{1} + 0\cdot n_{1}^{0}$.
\end{lemma}

\vspace{+1ex}
\noindent \textit{Proof}: If $g_{n}(m) = 0$ for some natural numbers $m$ and $n$ then:

\begin{eqnarray*}
If\ n > 1,\ then\ g_{n}(m)_{[(n+1)]} & = & 0\cdot(n+1)^{0} \\
If\ n > 2,\ then\ g_{(n-1)}(m)_{[n]} & = & 1\cdot(n)^{0} \\
If\ n > 3,\ then\ g_{(n-2)}(m)_{[(n-1)]} & = & 2\cdot(n-1)^{0} \\
If\ n > 5,\ then\ g_{(n-3)}(m)_{[(n-2)]} & = & 3\cdot(n-2)^{0} \\
\ldots \\
If\ n > (2r-1),\ then\ g_{(n-r)}(m)_{[(n-r+1)]} & = & r\cdot(n-r+1)^{0}.
\end{eqnarray*}

\vspace{+1ex}
\noindent Now, for some natural number $n_{1} > 1$, we must have either $n = (2n_{1}+1)$ or $n = 2n_{1}$. If $n = (2n_{1}+1)$ then:

\begin{eqnarray*}
g_{(n-n_{1})}(m)+1 & = & n_{1}\cdot(n-n_{1}+1)^{0}+1\cdot(n-n_{1}+1)^{0} \\
& = & (n_{1}+1)\cdot(n_{1}+2)^{0}
\end{eqnarray*}

\vspace{+1ex}
However, since $(n_{1}+1)\cdot(n_{1}+2)^{0}$ is not the result of any Goodstein bumping operation, we cannot have $n = (2n_{1}+1)$. 

\vspace{+1ex}
Hence, $n = 2n_{1}$ and:

\begin{eqnarray*}
g_{(n-n_{1})}(m)_{[(n-n_{1}+1)]} & = & n_{1}\cdot(n-n_{1}+1)^{0} \\
g_{n_{1}}(m)_{[(n_{1}+1)]} & = & n_{1}\cdot(n_{1}+1)^{0} \\
g_{(n_{1}-1)}(m)_{[n_{1}]} & = & 1\cdot n_{1}^{1} + 0\cdot n_{1}^{0}
\end{eqnarray*}

\vspace{+1ex}
\noindent The lemma follows. \hfill $\Box$

\vspace{+1ex}
We further have:

\begin{corollary} If the $n^{th}$ term $g_{n}(m)$ of the Goodstein sequence $G(m)$ is $0$, then $n = 2n_{1}$ and, for $n_{1} \leq k \leq n$, we have that $g_{k}(m) = (n-k)$.    
\end{corollary}

\subsection{The fifth Goodstein sequence lemma}
\label{lem:1.7}

\begin{lemma}
\label{sec:lem5}
If the $n^{th}$ term $g_{n}(m)$ of the Goodstein sequence $G(m)$ is $0$, then $n = 2(2n_{2}+1)$ for some natural number $n_{2}$ and $g_{(n_{2}-1)}(m)_{[n_{2}]} = 2\cdot n_{2}^{1} + 0\cdot n_{2}^{0}$.
\end{lemma}

\vspace{+1ex}
\noindent \textit{Proof}: Arguing as before, we have that:

\begin{eqnarray*}
If\ n_{1} > 2,\ then\ g_{(n_{1}-1)}(m)_{[n_{1}]} & = & 1\cdot n_{1}^{1} + 0\cdot n_{1}^{0} \\
If\ n_{1} > 3,\ then\ g_{(n_{1}-2)}(m)_{[(n_{1})-1]} & = & 1\cdot (n_{1}-1)^{1} + 0\cdot(n_{1}-1)^{0} \\
If\ n_{1} > 4,\ then\ g_{(n_{1}-3)}(m)_{[(n_{1})-2]} & = & 1\cdot (n_{1}-2)^{1} + 0\cdot(n_{1}-2)^{0} \\
If\ n_{1} > 6,\ then\ g_{(n_{1}-4)}(m)_{[(n_{1})-3]} & = & 1\cdot (n_{1}-3)^{1} + 0\cdot(n_{1}-3)^{0} \\
\ldots \\
If\ n_{1} > 2r,\ then\ g_{(n_{1}-r)}(m)_{[(n_{1})-r+1]} & = & 1\cdot (n_{1}-r+1)^{1} + 0\cdot(n_{1}-r+1)^{0}
\end{eqnarray*}

\vspace{+1ex}
We thus have that $(n_{1}-1)$ must also be even and, if $(n_{1}-1) = 2n_{2}$, then:

\begin{quote}
If $n_{2} > 3$ then $g_{(n_{2}-1)}(m)_{[n_{2}]} = 2\cdot n_{2}^{1} + 0\cdot n_{2}^{0}$.
\end{quote}

\noindent The lemma follows. \hfill $\Box$

\begin{corollary}
If the $n^{th}$ term $g_{n}(m)$ of the Goodstein sequence $G(m)$ is $0$, then $n = 2(2n_{2}+1)$ for some natural number $n_{2}$ and, for $n_{2} \leq k < 2n_{2}$, we have that $g_{k}(m) = g_{(k+1)}(m)$.
\end{corollary}

\subsection{The sixth Goodstein sequence lemma}
\label{lem:1.8}

\begin{lemma}
\label{sec:lem6}
$(p>m) \rightarrow (g_{n}(p)\geq g_{n}(m))$
\end{lemma}

\vspace{+1ex}
\noindent \textit{Proof}: We have that:

\begin{eqnarray*}
p>m & \equiv & g_{1}(p)>g_{n}(m) \\
g_{(n)}(p)_{[(n+1)]} > g_{(n)}(m)_{[(n+1)]} & \equiv & (\sum_{i=0}^{l_{p}} a_{i}.(n+1)^{i_{[(n+1)]}}) \\ & & > (\sum_{i=0}^{l_{m}} a_{i}.(n+1)^{i_{[(n+1)]}}) \\
& \rightarrow & (\sum_{i=0}^{l_{p}} a_{i}.(n+2)^{i_{[(n+1) \hookrightarrow (n+2)]}}) \\ & & > (\sum_{i=0}^{l_{m}} a_{i}.(n+2)^{i_{[(n+1) \hookrightarrow (n+2)]}}) \\
& \rightarrow & (\sum_{i=0}^{l_{p}} a_{i}.(n+2)^{i_{[(n+1) \hookrightarrow (n+2)]}}-1) \\ & & \geq (\sum_{i=0}^{l_{m}} a_{i}.(n+2)^{i_{[(n+1) \hookrightarrow (n+2)]}}-1) \\
& \rightarrow & g_{(n+1)}(p)_{[(n+2)]} \geq g_{(n+1)}(m)_{[(n+2)]}
\end{eqnarray*}

\noindent The lemma follows by induction. \hfill $\Box$

\section{Three Goodstein sequence theorems}
\label{seqtheo:1}

It immediately follows from the above lemmas that:

\begin{theorem}
\label{theo:1}
(First Goodstein sequence theorem): If the nth term $g_{n}(m)$ of the Goodstein sequence $G(m)$ is $0$, then $n = 2(2n_{2}+1)$ for some natural number $n_{2}$ and:

\begin{eqnarray*}
g_{k}(m) & = & (n-k)\ for\ (2n_{2}+1) \leq k \leq n \\
g_{k}(m) & = & g_{(k+1)}(m)\ for\ n_{2} \leq k < 2n_{2} \\
g_{k}(m) & < & g_{(k+1)}(m)\ for\ 1 \leq k < n_{2} 
\end{eqnarray*}
\end{theorem}

Theorem \ref{theo:1} highlights the characteristic structure of all terminating Goodstein sequences as schematically depicted below for $G(p)$ and $G(m)$ where $p>m>4$ \footnote{The characteristic structure of all terminating Goodstein sequences, expressed by this theorem, was investigated interestingly by R.E.S. \cite{Re03}, and visualised intriguingly in a striking computer-generated plot of $G(m)$ for $m \geq 4$ as mentioned \href{http://alixcomsi.com/RES_Goodstein_Sequences.htm}{here}.}.

\vspace{+2ex}
\begin{picture}(300,200)
\put(0,0){\vector(1,0){320}}
\put(0,0){\vector(0,1){200}}

\put(10,60){\circle*{3}}
\qbezier(10,60)(20,170)(80,190)

\put(80,190){\circle*{3}}
\put(80,190){\line(70,0){70}}

\put(150,190){\circle*{3}}
\qbezier(150,190)(290,0)(290,0)
\put(290,0){\circle*{3}}

\put(10,10){\circle*{3}}
\qbezier(10,10)(25,120)(60,130)

\put(60,130){\circle*{3}}
\put(60,130){\line(50,0){50}}

\put(110,130){\circle*{3}}
\qbezier(110,130)(210,0)(210,0)
\put(210,0){\circle*{3}}

\put(3,55){\scriptsize{$g_{1}(p)=p$}}
\put(60,195){\scriptsize{$g_{s}(p)= \ldots$}}
\put(120,195){\scriptsize{$\ldots = g_{(2s-1)}(p)$}}
\put(270,-6){\scriptsize{$g_{q}(p)=0$}}

\put(3,5){\scriptsize{$g_{1}(m)=m$}}
\put(45,135){\scriptsize{$g_{r}(m)= \ldots$}}
\put(90,135){\scriptsize{$\ldots = g_{(2r-1)}(m)$}}
\put(190,-6){\scriptsize{$g_{n}(m)=0$}}

\put(-21,1){\scriptsize{Value}}
\put(-4,8){\vector(0,1){15}}

\put(1,-6){\scriptsize{Step}}
\put(18,-4){\vector(1,0){15}}

\put(80,-6){\tiny{(Not to scale)}}

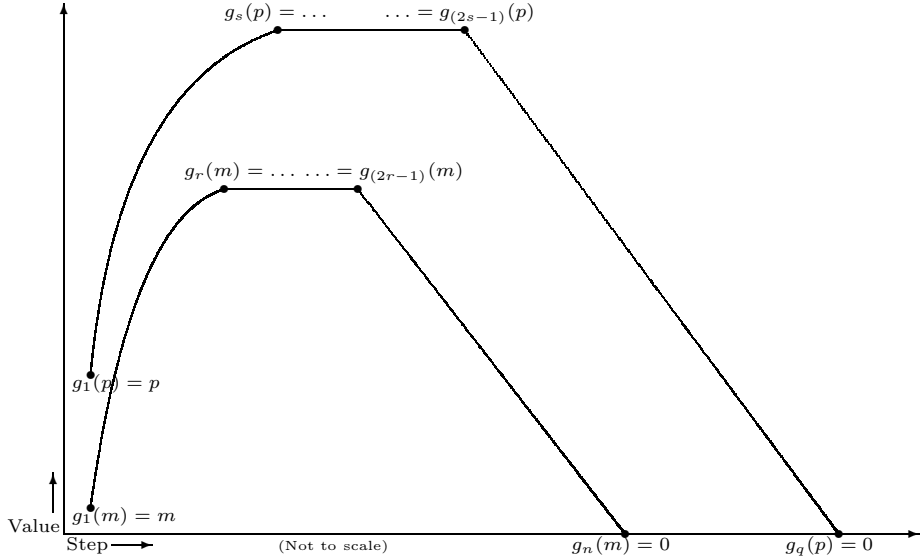
\put(0,-25){\footnotesize{Fig 1: Characteristic structure of terminating Goodstein Sequences: $m<p$}}
\end{picture}
                                                    
\vspace{+8ex}                               
\noindent We can further conclude that:                                                               

\begin{theorem}     
Any sequence generated by the iterated application of the complete Goodstein operation on the hereditary representation of any natural number in any natural number base cannot oscillate.
\end{theorem}

\begin{corollary}
For any given natural number $m$, the Goodstein sequence $G(m)$ of natural numbers converges if, and only if, it is bounded finitely in Peano Arithmetic.
\end{corollary}

We note next that if $H(k)$ is the sequence of finite natural numbers generated by iterated application of the complete Goodstein operation on the hereditary representation of $k^{k}$ in the base $[k]$, then:

\begin{lemma}
\label{lem:theo3}
For any natural number $k$, the sequence $H(k) = \{h_{1}(k), h_{2}(k), h_{3}(k), \ldots\}$ defined by:

\begin{eqnarray*}
h_{1}(k)_{[k]} & = & k^{k_{[k]}} \\
h_{n+1}(k)_{[(k+n)]} & = & h_{n}(k)_{[(k+n-1) \hookrightarrow\ (k+n)]} - 1
\end{eqnarray*}

\vspace{+1ex}
\noindent is a finite sequence that terminates with $0$.
\end{lemma}

\vspace{+1ex}
\noindent \textit{Proof}: By definition:

\begin{eqnarray*}
h_{1}(k)_{[k]} & = & 1\cdot k^{k_{[k]}} + 0\cdot k^{(k-1)_{[k]}} + 0\cdot k^{(k-2)_{[k]}} + \ldots + 0\cdot k^{0_{[k]}} \\
h_{2}(k)_{[k+1]} & = & (1\cdot (k+1)^{k_{[k \hookrightarrow\ (k+1)]}} + 0\cdot (k+1)^{(k-1)_{[k \hookrightarrow\ (k+1)]}} + \\
& & 0\cdot (k+1)^{(k-2)_{[k \hookrightarrow\ (k+1)]}} + \ldots + 0\cdot (k+1)^{0_{[k \hookrightarrow\ (k+1)]}}) - 1 \\
& = & (1\cdot (k+1)^{(k+1)_{[(k+1)]}} + 0\cdot (k+1)^{(k-1)_{[(k+1)]}} + \\
& & 0\cdot (k+1)^{(k-2)_{[(k+1)]}} + \ldots + 0\cdot (k+1)^{0_{[(k+1)]}}) - 1 \\
& = & k\cdot (k+1)^{k_{[(k+1)]}} + k\cdot (k+1)^{(k-1)_{[(k+1)]}} + \\
& & k\cdot (k+1)^{(k-2)_{[(k+1)]}} + \ldots + k\cdot (k+1)^{0_{[(k+1)]}} \\
\ldots
\end{eqnarray*}

\vspace{+1ex}
Thus, $H(k)$ is a sequence of natural numbers each of whose members is such that the rank of the hereditary representation of any successor-member is less than the base of the representation; in other words, the complete Goodstein operation leaves the exponents in the hereditary representation of any successor-member unchanged (as in Goodstein's argument with transfinite ordinals).

Now, we note that a complete Goodstein operation on any number of the form:

\begin{eqnarray*}
\sum_{i=0}^{l} a_{i}k^{i}
\end{eqnarray*}

\vspace{+1ex}
\noindent where $1 \leq a_{l}, l < k$ and $a_{i} = 0$ for $0 \leq i < l$:

\begin{quote}
(i)	either yields, if $a_{l} > 1$, a number of the same rank $l$ but with a reduced co-efficient of the highest power of the base, such as:

\begin{eqnarray*}
\sum_{i=0}^{l} c_{i}(k+1)^{i}
\end{eqnarray*}

\vspace{+1ex}
\noindent where $c_{l} = (a_{l}-1)$ and $0 < c_{i} < (k+1)$ for $0 \leq i < l$; 

\vspace{+1ex}
(ii)	or yields, if $a_{l} = 1$, a number of reduced rank such as:

\begin{eqnarray*}
\sum_{i=0}^{l-1} c_{i}(k+1)^{i}
\end{eqnarray*}

\vspace{+1ex}
\noindent where $0 < c_{i} < (k+1)$ for $0 \leq i < (l-1)$.
\end{quote}

\vspace{+1ex}
In either case, it can be shown by finite induction that iterated application of the complete Goodstein operation must eventually reduce the rank of some member of the sequence to $0$. By Lemma \ref{lem:1.3} the sequence must therefore terminate with $0$.

The lemma follows. \hfill $\Box$

\vspace{+1ex}
We, thus, have that:

\begin{theorem}
(Third Goodstein sequence theorem): For any given natural number $m$ the Goodstein sequence $G(m)$ converges finitely if, and only if, there is some natural number $k$ such that:

\begin{eqnarray*}
g_{n}(m) > 0 \rightarrow (h_{n}(k) > g_{n}(m))
\end{eqnarray*}
\end{theorem}

\vspace{+1ex}
\noindent \textit{Proof}: We note firstly that, if there is a convergent sequence $H(k)$ as defined above, then $G(m)$ is a finitely bounded sequence and so, by the preceding lemmas, it must converge to $0$. 

Secondly, if $G(m)$ converges finitely in $n$ terms, and the largest member of the sequence is $p$, we can take $k = max(n, p)$. By lemma \ref{lem:theo3}, $h_{2q}(k) = 0$ for some $q$. 

Since, by the first Goodstein sequence theorem \ref{theo:1}, the largest term in $H(k)$ is $h_{q}(k) = q$ and we also have that $q > k^{k}$, it follows that $q > n$. Hence for any $n$ such that $g_{n}(m) \neq 0$ we have that $h_{n}(k) > g_{n}(m)$. \hfill $\Box$

\section{Goodstein's implicit Thesis}
\label{thesis:1}

We can now express Goodstein's implicit thesis as: 

\begin{thesis}
(The Goodstein Thesis): For any given natural numbers $m$ and $k$, and partial sequences $\{g_{n}(m): 1 \leq n \leq k\}$, i.e.:

\begin{eqnarray*}
\{g_{1}(m)_{[2]},\ g_{2}(m)_{[3]},\ \ldots,\ g_{k}(m)_{[(k+1)]}\}
\end{eqnarray*}

\vspace{+1ex}
\noindent and $\{g_{n}(m)_{[(n+1)\ \hookrightarrow\ u_{k}]}: 1 \leq n \leq k\}$, i.e.:

\begin{eqnarray*}
\{g_{1}(m)_{[2\ \hookrightarrow\ u_{k}]},\ g_{2}(m)_{[3\ \hookrightarrow\ u_{k}]},\ \ldots,\ g_{k}(m)_{[(k+1)\ \hookrightarrow\ u_{k}]}\}
\end{eqnarray*}

\vspace{+1ex}
\noindent where:

\begin{quote}
(i)	$[u_{k}]$ is the base of the largest term in the partial sequence $\{g_{n}(m): 1 \leq n \leq k\}$; 

(ii) and $g_{n}(m)_{[(n+1)\ \hookrightarrow\ u_{k}]}$ is the natural number obtained from the hereditary representation of $g_{n}(m)_{[(n+1))]}$ by syntactically replacing the base $[(n+1)]$ by the base $[u_{k}]$;
\end{quote}

we have that:

\begin{eqnarray*}
(Lt\ u_{k}) < \omega
\end{eqnarray*}
\end{thesis}

\vspace{+1ex}
Now, by the first Goodstein sequence theorem there is no natural number base $[u_{k}]$ for which we have that, for all $n \geq 1$:

\begin{eqnarray*}
g_{n}(m)_{[(n+1)\ \hookrightarrow\ u_{k}]} < g_{(n+1)}(m)_{[(n+2)\ \hookrightarrow\ u_{k}]} 
\end{eqnarray*}

\vspace{+1ex}
Thus Goodstein's argument appeals to properties of sequences of transfinite ordinals in Cantor Arithmetic that do not correspond to any arithmetical properties of their corresponding sequences of natural numbers in Peano Arithmetic. 

It follows we cannot prima facie conclude that, by simply replacing a constructive natural number base $[u_{k}]$ by the non-constructive ordinal base $[\omega]$, Goodstein's argument establishes Goodstein's Theorem as a true assertion of Peano Arithmetic under its standard interpretation in a constructive, and intuitionistically unobjectionable, way. 

Such a conclusion must therefore implicitly appeal to some non-constructive, and counter-intuitive, assumption that needs to be expressed explicitly.

\section{A Goodstein functional sequence}
\label{func:1}

The above point is illustrated better if we define the sequence of functions:

\begin{eqnarray*}
\{g_{1}(m)_{[2\ \hookrightarrow\ x]},\ g_{2}(m)_{[3\ \hookrightarrow\ x]},\ \ldots,\ g_{k}(m)_{[(k+1)\ \hookrightarrow\ x]}\}
\end{eqnarray*}

\vspace{+1ex}
\noindent as the Goodstein functional sequence of $m$.

\vspace{+1ex}
Now, if $g_{n}(m)_{[(n+1)\ \hookrightarrow\ x]}$ has a constant term $c$ where $x > c > 0$, we have:

\begin{eqnarray*}
g_{n}(m)_{[(n+1)\ \hookrightarrow\ x]} - g_{(n+1)}(m)_{[(n+2)\ \hookrightarrow\ x]} = 1
\end{eqnarray*}

\vspace{+1ex}
\noindent whilst if the lowest power of $x$ in $g_{n}(m)_{[(n+1)\ \hookrightarrow\ x]}$ is $c \cdot x^{a}$, where $c$ and $a$ are constants such that $a > 0$ and $x > c > 0$, then:

\begin{eqnarray*}
g_{n}(m)_{[(n+1)\ \hookrightarrow\ x]} - g_{(n+1)}(m)_{[(n+2)\ \hookrightarrow\ x]} \geq (x-2)x^{(a-1)}
\end{eqnarray*}

\vspace{+1ex}
Hence, for all natural numbers $u > 2$:

\begin{eqnarray*}
\{g_{n}(m)_{[(n+1)\ \hookrightarrow\ u]}\ :\ n \geq 1\}
\end{eqnarray*}

\vspace{+1ex}
\noindent is a finite, decreasing, sequence of natural numbers.

Now, if the Goodstein sequence $G(m)$ terminates after $l$ terms then, for any $n \geq 1$:

\begin{eqnarray*}
g_{n}(m)_{[(n+1)\ \hookrightarrow\ l]}\ \geq g_{n}(m)_{[(n+1)]}
\end{eqnarray*}

\vspace{+1ex}
However, if $G(m)$ does not terminate then, for any given $u > 2$, there can be no decreasing sequence of natural numbers $\{g_{n}(m)_{[(n+1)\ \hookrightarrow\ u]}\ :\ n \geq 1\}$ such that $g_{n}(m)_{[(n+1)\ \hookrightarrow\ l]}\ \geq g_{n}(m)_{[(n+1)]}$ for any $n \geq 1$.

However in either case---treating $g_{n}(m)_{[(n_{o}+1)\ \hookrightarrow\ x]}$ as a function over the ordinals---there is always a decreasing sequence of transfinite ordinals:

\begin{eqnarray*}
\{g_{n}(m_{o})_{[(n_{o}+1)\ \hookrightarrow\ \omega]}\ :\ n \geq 1\}
\end{eqnarray*}

\vspace{+1ex}
\noindent where $m_{o},\ n_{o}$ are the finite ordinals corresponding to the natural number $m,\ n$.

\section{Formal mathematical objects}
\label{obj:1}

Now, as noted above, Goodstein's argument does not directly address this issue of the upper bound of a Goodstein sequence in Peano Arithmetic. Instead it indirectly---and implicitly---concludes the existence of such a bound since there is a finite, recursively well-defined, decreasing sequence $\{g_{n}(m_{o})_{[(n_{o}+1_{o})\ \hookrightarrow\ \omega]}\}$ of transfinite ordinals in Cantor Arithmetic such that, for any $n$ such that $g_{n}(m_{o})_{[(n_{o}+1_{o})]} \neq 0_{o}$, we have that $g_{n}(m_{o})_{[(n_{o}+1_{o})\ \hookrightarrow\ \omega]} \geq g_{n}(m_{o})_{[(n_{o}+1_{o})]}$.

However, in the absence of a constructive proof---or meta-proof\footnote{We note that such a meta-proof need not be provable in PA; in other words, the truth---or falsity---of Goodstein's Theorem in the standard interpretation of PA does not depend on the existence of a PA proof sequence for the Theorem (or for its negation).}---to the contrary, we must admit the possibility that some Goodstein sequence of natural numbers in Peano Arithmetic does not converge. The latter would be the case if Goodstein's assumption---that we can recursively define an ordinal sequence $\{g_{n}(m_{o})_{[(n_{o}+1_{o})\ \hookrightarrow\ \omega]}\}$ as a well-defined, formal, mathematical object in Cantor Arithmetic---is invalid. Consequently, Goodstein's argument---that this ordinal sequence is a well-defined, finite\footnote{As ([Me64 previously noted, for the sequence to be termed as finite, it must be a well-defined set in Cantor Arithmetic], p184).}, decreasing, set in any putative model of Cantor Arithmetic---would be vacuously true. 

It follows that, if we admit such a possibility, then we cannot treat the standard interpretations of Goodstein's argument as establishing Goodstein's Theorem definitively.

We can express the above reservation more precisely. Goodstein's argument ignores the possibility that the recursively defined ordinal sequence:

\begin{eqnarray*}
\{g_{n}(m_{o})_{[(n_{o}+1_{o})\ \hookrightarrow\ \omega]}\}
\end{eqnarray*}

\vspace{+1ex}
\noindent may not be a formal mathematical object in Cantor Arithmetic (i.e., in a set theory such as, say, ZFC), in the following sense:

\begin{quote}
\begin{definition}
A primitive mathematical object is any symbol for an individual constant, predicate letter, or a function letter\footnote{cf. \cite{Me64}, p46; also p1, p10.} which is defined as a primitive symbol of a formal mathematical language.
\end{definition}

\begin{definition}
A formal mathematical object is any symbol for an individual constant, predicate letter, or a function letter that is either a primitive mathematical object, or that can be introduced through definition\footnote{cf. \cite{Me64}, p82.} into a formal mathematical language without inviting inconsistency\footnote{We take Mendelson's Corollary 1.15 (\cite{Me64}, p37) as the classical meta-definition of consistency.}.
\end{definition}

\begin{definition}
A mathematical object is any symbol that is either a primitive mathematical object or a formal mathematical object.
\end{definition}

\begin{definition}
A set is the range of any function whose function letter is a formal mathematical object.
\end{definition}
\end{quote}

Consideration of formal mathematical objects in more detail would lie outside the immediate scope of this investigation, whose limited aim is to establish that, prima facie, there are sufficient grounds for arguing that the standard interpretations of Goodstein's argument ought not to be accepted as definitive. 

Nevertheless we note that in \cite{An02b} we consider the existence of a primitive recursive number-theoretic relation that is intuitively decidable constructively, but which cannot be introduced through definition as a formal mathematical object into the formal system of Peano Arithmetic PA without inviting inconsistency. Ipso facto, such a relation cannot be introduced through definition as a formal mathematical object into any Axiomatic Set Theory, such as ZFC, in which the axioms of PA interpret as theorems. Hence it is not a formal mathematical object and the range of its characteristic function is not a recursively enumerable set.

Since recursive number-theoretic functions and relations are classically accepted as amongst the most basic building blocks for defining constructive, and intuitionistically unobjectionable, mathematical objects, we cannot prima facie accept Goodstein's argument as sufficient for establishing the recursively defined---and admittedly non-constructive---ordinal sequence $\{g_{n}(m_{o})_{[(n_{o}+1_{o})\ \hookrightarrow\ \omega]}\}$ as a formal mathematical object in Cantor Arithmetic.

\section{Conclusion}
\label{conc:1}

Goodstein's argument implicitly assumes that the recursively defined ordinal sequence $\{g_{n}(m_{o})_{[(n_{o}+1_{o})\ \hookrightarrow\ \omega]}\}$ is a formal mathematical object in Cantor Arithmetic. In other words, it implicitly assumes the existence of a well-defined set of transfinite ordinals in Cantor Arithmetic which has properties corresponding to the properties required of the number-theoretic sequence $H(k)$ that is defined in Lemma \ref{lem:theo3} above.

Since, as argued in \cite{An02b}, such an assumption need not necessarily hold, Goodstein's Theorem can reasonably be viewed as a number-theoretic proposition whose truth in the standard interpretation of any formal system of Peano Arithmetic has simply been asserted as a non-verifiable consequence of a non-constructive argument.

We conclude that in the absence of a constructive and intuitionistically unobjectionable proof---or meta-proof---that any given Goodstein sequence is bounded in Peano Arithmetic, the standard interpretation of Goodstein's Theorem as a number-theoretic assertion that is consistent with any formal system of Peano Arithmetic, such as standard PA, ought not to be accepted as definitive.

\noindent \tiny{Authors postal address: 32 Agarwal House, D Road, Churchgate, Mumbai - 400 020, Maharashtra, India.\ Email: re@alixcomsi.com, anandb@vsnl.com.}

\vspace{+1ex}
(Updated: Friday 5th August 2004 11:12:03 PM IST by re@alixcomsi.com)

\end{document}